%Two Notes on Notation
\magnification\magstephalf
\baselineskip 14pt
\parskip 2pt

\input picmac

\font\eusm=eusm10
\font\eufr=eufm10
\font\sc=cmcsc10 %use lower case as (Monthly)

\def\TeX{T\hbox{\hskip-.1667em\lower.424ex\hbox{E}\hskip-.125em X}}

\def\Piit{{\mit\Pi}}

\def\pfbox
  {\hbox{\hskip 3pt
\vbox{\hrule\hbox to 7pt{\vrule height 7pt\hfill\vrule}
  \hrule}}\hskip3pt}

\def\bib{\par\noindent\hangindent 20pt}

\centerline{\bf Two Notes on Notation}
\centerline{by Donald E. Knuth\footnote{}{%
\kern-\parindent
This research
was supported in part by National Science Foundation grant
CCR-86-10181.}}
\centerline{\sl Computer Science Department, Stanford University}

\bigskip
Mathematical notation evolves like all languages do. As new
experiments are made, we sometimes witness the survival of the
fittest, sometimes the survival of the most familiar. A~healthy
conservatism keeps things from changing too rapidly; a~healthy
radicalism keeps things in tune with new theoretical emphases. Our
mathematical language continues to improve, just as ``the $d$-ism of
Leibniz overtook the dotage of Newton'' in past centuries
[4, Chapter~4].

In 1970 I began teaching a class at Stanford University entitled
Concrete Mathematics. The students and~I studied how to manipulate
formulas in continuous and discrete mathematics,  and the problems we
investigated were often inspired by new developments in computer
science. As the years went by we 
began to see that a few changes in notational traditions would
greatly facilitate our work. The notes from that class have recently
been published in a book~[15], and as I~wrote the final drafts of
that book I~learned to my surprise that two of the notations we had
been using were considerably more useful than I~had previously
realized. The ideas ``clicked'' so well, in fact, that I've decided to
write this article, blatantly attempting to promote these notations
among the mathematicians who have no use for~[15]. I~hope that within
five years everybody will be able to use these notations in published
papers without needing to explain what they mean.

The notations I'm talking about are (1) Iverson's convention for
characteristic functions; and (2)~the ``right'' notation for Stirling
numbers, at last. 

\bigskip\noindent
{\bf 1. Iverson's convention}.\enspace
The first notational development I~want to discuss was introduced by
Kenneth E.~Iverson in the early~60s, on page~11 of the pioneering 
book~[21] 
that led to his well known {\sl APL}.

{\narrower\smallskip\noindent
\llap{``}If $\alpha$ and $\beta$ are arbitrary entities and {\eusm R}
is any relation defined on them, the {\it relational statement\/}
$(a\hbox{\eusm R}b)$ is a logical variable which is true (equal to~1)
if and only if $\alpha$~stands in the relation~{\eusm R} to~$\beta$.
For example, if $x$~is any real number, then the function
$$(x>0) \;-\;(x<0)$$ 
(commonly called the {\it sign function\/} or sgn~$x$) assumes the
values 1, 0, or~$-1$ according as $x$~is strictly positive, 0, or
strictly negative.''
\smallskip}

When I read that, long ago, I~found it mildly interesting but not
especially significant. I~began using his convention informally but
infrequently, in class discussions and in private notes. 
I~allowed it to slip, undefined, into an obscure corner of
 one of my books (see page~117 of~[16]). But when I~prepared the final
manuscript of~[15], I~began to notice that Iverson's idea led to
substantial improvements in exposition and in technique.

Before I can explain why the notation now works so well for me, I~need
to say a few words about the manipulation of sums and summands.
I~realized long ago that ``boundary conditions'' on indices of
summation are often a handicap and a waste of time. Instead of writing
$$(1+z)^n=\sum_{k=0}^n{n\choose k}\,z^k\,,\eqno(1.1)$$
it is much better to write
$$(1+z)^n=\sum_k{n\choose k}\,z^k\,;\eqno(1.2)$$
the sum now extends over all integers~$k$, but only
 finitely many terms are nonzero. The
second formula (1.2) is instantly converted to other forms:
$$(1+z)^n=\sum_k{n\choose k}\,z^k=\sum_k{n\choose k+1}\,z^{k+1}=\sum_k
{n\choose \lfloor n/2\rfloor -k}\,z^{\lfloor n/2\rfloor
-k}\,;\eqno(1.3)$$
by contrast, we must work harder when dealing with (1.1), because we
have to think about the limits:
$$(1+z)^n=\sum_{k=0}^n{n\choose k}\,z^k=\sum_{k=-1}^{n-1}{n\choose
k+1}\,z^{k+1}=\sum_{k=-\lceil n/2\rceil}^{\lfloor n/2\rfloor}
{n\choose \rfloor n/2\rfloor -k}\,z^{\lfloor n/2\rfloor
-k}\,.\eqno(1.4)$$
Furthermore, (1.2) and (1.3) make sense also when $n$ is not a
positive integer.

Even when limits are necessary, it is best to keep them as simple as
possible. For example, it's almost always a mistake to write
$$\sum_{k=2}^{n-1}k(k-1)(n-k)\qquad\hbox{instead of}\qquad
\sum_{k=0}^n k(k-1)(n-k)\,;\eqno(1.5)$$
the additional zero terms are more helpful than harmful (and the
former sum is problematical when $n=0$, 1, or~2).

Finally it dawned on me that Iverson's convention allows us to write
{\it any\/} sum as an infinite sum without limits: If $P(k)$ is any
property of the integer~$k$, we have
$$\sum_{P(k)}f(k)=\sum_kf(k)\,[P(k)]\,.\eqno(1.6)$$
For example, the sums in (1.5) become
$$\sum_k k(k-1)(n-k)\;[0\leq k\leq n]=\sum_k k(k-1)(n-k)\;[k\geq 0]\,
[k\leq n]\,.\eqno(1.7)$$
(At the time I made this observation, I had forgotten that Iverson originally
defined his convention only for single relational operators enclosed
in parentheses; I~began to put {\it arbitrary\/} logical statements in
square brackets, and to assume that this would produce the value~0
or~1.)
In this particular case nothing much has been gained when passing from
(1.5) to (1.7), although we might be able to make use of identities
like
$$k\;[k\geq 0]\;=\;k\;[k\geq 1]\,.\eqno(1.8)$$
But in general, the ability to manipulate ``on the line'' instead of
``below the line'' turns out to be a great advantage.

For example, in my first book~[25] I~had found it necessary to
include the rule
$$\sum_{k\in A}f(k)+\sum_{k\in B}f(k)\ =\sum_{k\in A\cup
B}f(k)\ +\sum_{k\in A\cap B}f(k)\eqno(1.9)$$
as a separate axiom for $\sum$ manipulation. But this axiom is
unnecessary in~[15], because it can be derived easily from other
basic laws: The left-hand side is
$$\eqalign{\sum_{k\in A}f(k)+\sum_{k\in B}f(k)&=\sum_kf(k)\,[k\in A]
+\sum_kf(k)\,[k\in B]\cr
\noalign{\smallskip}
&=\sum_kf(k)\,([k\in A]+[k\in B])\cr}$$
and the right-hand side is the same, because we have
$$[k\in A]+[k\in B]=[k\in A\cup B]+[k\in A\cap B]\,.\eqno(1.10)$$

The interchange of summation order in multiple sums also comes out simpler
now. I~used to have trouble understanding and/or explaining why
$$\sum_{j=1}^n\,\sum_{k=1}^jf(j,k)=\sum_{k=1}^n\,\sum_{j=k}^nf(j,k)\,;
\eqno(1.11)$$
but now it's easy for me to see that the left-hand sum is
$$\eqalign{\sum_{j,k}f(j,k)\,[1\leq j\leq n]\,[1\leq k\leq j]
&=\sum_{j,k}f(j,k)\,[1\leq k\leq j\leq n]\cr
\noalign{\smallskip}
&=\sum_{j,k}f(j,k)\,[1\leq k\leq n]\,[k\leq j\leq n]\,,\cr}$$
and this is the right-hand sum.

Here's another example: We have
$$[k\;{\rm even}]=\sum_m\,[k=2m]\qquad{\rm and}\qquad [k\;{\rm odd}]
=\sum_m\,[k=2m+1]\,;\eqno(1.12)$$
therefore
$$\eqalignno{\sum_kf(k)&=\sum_kf(k)\,([k\;{\rm even}]+[k\;{\rm odd}])\cr
\noalign{\smallskip}
&=\sum_{k,m}f(k)\,[k=2m]+\sum_{k,m}f(k)\,[k=2m+1]\cr
\noalign{\smallskip}
&=\sum_mf(2m)+\sum_mf(2m+1)\,.&(1.13)\cr}$$
The result in (1.13) is hardly surprising; but I like to have
mechanical operations like this available so that I~can do
manipulations reliably, without thinking. Then I'm less apt to make
mistakes.

Let $\lg$ stand for logarithms to base~2. Then we have
$$\eqalignno{%
\sum_{k\geq 1}{n\choose\lfloor \lg k\rfloor}&=\sum_{k\geq 1}\,\sum_m
{n\choose m}\,\bigl[m=\lfloor\lg k\rfloor\bigr]\cr
\noalign{\smallskip}
&=\sum_{k,m}{n\choose m}\,[m\leq \lg k<m+1]\,[k\geq 1]\cr
\noalign{\smallskip}
&=\sum_{m,k}{n\choose m}\,[2^m\leq k<2^{m+1}]\,[k\geq 1]\cr
\noalign{\smallskip}
&=\sum_m{n\choose m}(2^{m+1}-2^m)\,[m\geq 0]\cr
\noalign{\smallskip}
&=\sum_m{n\choose m}2^m=3^n\,.&(1.14)\cr}$$

If we are doing infinite products we can use Iversonian brackets as
exponents:
$$\prod_{P(k)}f(k)=\prod_kf(k)^{[P(k)]}\,.\eqno(1.15)$$
For example, the largest squarefree divisor of $n$ is
$$\prod_p p^{\,[p\;{\rm prime}]\,[p\;{\rm divides}\;n]}\,.$$

Everybody is familiar with one special case of an Iverson-like
convention, the ``Kronecker delta'' symbol
$$\delta_{ik}=\cases{1\,,&$i=k$;\cr
\noalign{\smallskip}
0\,,&$i\not= k$.\cr}\eqno(1.16)$$
Leopold Kronecker introduced this notation in his work on bilinear 
forms~[30,
page 276] and in his lectures on determinants (see~[31, page 316]);
it soon became widespread. Many of his followers wrote~$\delta_j^k$,
which is a bit more ambiguous because it conflicts with ordinary
exponentiation. I~now prefer to write $[j=k]$ instead
of~$\delta_{jk}$,
because Iverson's convention is much more general. Although `$[j=k]$'
involves five written characters instead of the three
in~`$\delta_{jk}$', we lose nothing in common cases when
`$[j=k+1]$' takes the place of~`$\delta_{j(k+1)}$'.

Another familiar example of a 0--1 function, this time from continuous
mathematics, is Oliver Heaviside's unit step function $[x\geq 0]$. (See~[44]
and~[37] for expositions of Heaviside's methods.)
It is clear that Iverson's convention will be as useful with
integration as it is with summation, perhaps even more~so. I~have not
yet explored this in detail, because [15] deals mostly with sums.

It's interesting to look back into the history of mathematics and see
how there was a craving for such notations before they existed. For
example, an Italian count named Guglielmo Libri published several
papers in the 1830s concerning properties of the function~$0^{0^x}$. He
noted~[32] that~$0^x$ is either~0 (if $x>0$) or~1 (if $x=0$) or~$\infty$
(if $x<0$), hence
$$0^{0^x}=[x>0]\,.\eqno(1.17)$$
But of course he didn't have Iverson's convention to work with; he 
was pleased to discover a way to denote the discontinuous function
$[x>0]$ without leaving the realm of operations acceptable in his day.
He believed that ``la fonction $0^{0^{x-n}}$ est d'un grand usage dans
l'analyse math\'ematique.'' And he noted in [33] that his formulas ``ne
renferment aucune notation nouvelle. \dots Les formules qu'on obtient
de cette mani\`ere sont tr\`es simples, et rentrent dans l'alg\`ebre
ordinaire.''

Libri wrote, for example,
$$(1-0^{0^{-x}})(1-0^{0^{x-a}})$$
for the function $[0\leq x\leq a]$, and he gave the integral formula
$${2\over\pi}\int_0^{\infty}{dq\cos qx\over 1+q^2}=e^x\cdot
0^{0^{-x}}+e^{-x}
\bigl(1-0^{0^{-x}}\bigr)={e^x\over 0^{-x}+1}+{e^{-x}\over 0^x+1}\,.$$
(Of course, we would now write the value of that integral as
$e^{-\vert x\vert}$, but a simple notation for absolute value wasn't
introduced until many years later. I~believe that the first appearance
of~`$\vert z\vert$' for absolute value in Crelle's journal---the
journal containing Libri's papers~[32] and~[33]---occurred on page~227
of~[56] in 1881. Karl Weierstrass was the inventor of this notation,
which was applied at first only to complex numbers; Weierstrass seems
to have published it first in 1876~[55].)

Libri applied his $0^{0^x}$ function to number theory by exhibiting a
complicated way to describe the fact that $x$ is a divisor of~$m$. In
essence, he gave the following recursive formulation: Let 
$P_0(x)=1$ and for $k>0$ let
$$P_k(x)=0^{0^{x-k}}P_0(x)-0^{0^{x-k+1}}P_1(x)-\cdots
-0^{0^{x-1}}P_{k-1}(x)\,.$$ 
Then the quantity
$${1-m\cdot 0^{0^{x-m}}P_0(x)-(m-1)\,0^{0^{x-m+1}}P_1(x)-\cdots -
2\cdot 0^{0^{x-2}}P_{m-2}(x)-0^{0^{x-1}}P_{m-1}(x)\over x}$$
turns out to equal 1 if $x$ divides~$m$, otherwise it is~0. (One way
to prove this, Iverson-wise, is to 
replace $0^{0^{x-k}}$ in Libri's formulas by $[x>k]$, and to show
first by induction that  $P_k(x)=[x\;{\rm
divides}\;k]-[x\;{\rm divides}\;k-1]$ for all $k>0$. Then if
$a_k(x)=k\,[x>k]$, we have
$$\eqalign{\sum_{k=0}^{m-1}a_{m-k}(x)P_k(x)
&=\sum_{k=0}^{m-1}a_{m-k}(x)\,([x\;{\rm divides}\; k]-[x\;{\rm
divides}\;k-1])\cr 
\noalign{\smallskip}
&=\sum_{k=0}^{m-1}\,[x\;{\rm divides}\;
k]\,\bigl(a_{m-k}(x)-a_{m-k-1}(x)\bigr)\,.\cr}$$
If the positive integer $x$ is not a divisor of $m$, the terms of this
new sum are zero except when $m-k=m\bmod x$, when we have
$a_{m-k}(x)-a_{m-k-1}(x)=1$. On the other hand if $x$~is a divisor
of~$m$, the only nonvanishing term occurs for $m-k=x$, when we have
$a_{m-k}(x)-a_{m-k-1}(x)=0-(x-1)$. Hence the sum is $1-x\,[x\;{\rm
divides}\;m]$. Libri obtained his complicated formula by a less direct method,
applying Newton's identities to compute the sum of the $m$th powers of
the roots of the equation $t^{x-1}+t^{x-2}+\cdots +1=0$.)

Evidently Libri's main purpose was to show that unlikely functions can
be expressed in algebraic terms, somewhat as we might wish to show
that some complex functions can be computed by a Turing Machine.
``Give me the function~$0^{0^x}$, and I'll give you an expression for
$[x\;{\rm divides}\;m]$.'' But our goal with Iverson's notation is, by
contrast, to find a simple and natural way to express quantities that
help us solve problems. If we need a function that is~1 if and only if
$x$~divides~$m$, we can now write $[x\;{\rm divides}\;m]$.

Some of Libri's papers are still well remembered, but [32] and [33]
are not. I~found no mention of them in {\sl Science Citation Index},
after searching through all years of that index available in our
library (1955 to date). However, the paper~[33]
 did produce several ripples in
mathematical waters when it originally appeared, because it stirred up
a controversy about whether $0^0$ is defined. Most mathematicians agreed
that $0^0=1$, but Cauchy~[5, page~70] had listed~$0^0$ together with
other expressions like $0/0$ and $\infty -\infty$ in a table of
undefined forms. Libri's justification for the equation $0^0=1$ was
far from convincing, and a commentator who signed his name
simply~``S'' rose to the attack~[45]. August M\"obius~[36] defended
Libri, by presenting his former professor's reason for believing that
$0^0=1$ (basically a proof that $\lim_{x\rightarrow 0+}x^x=1$).
M\"obius also went further and presented a supposed proof that
$\lim_{x\rightarrow 0+}f(x)^{g(x)}=1$ whenever 
\hbox{$\lim_{x\rightarrow
0+}f(x)=\lim_{x\rightarrow 0+}g(x)=0$}. Of course ``S''~then asked~[3]
whether M\"obius knew about functions such as $f(x)=e^{-1/x}$ and
$g(x)=x$. (And paper~[36] was quietly omitted from the historical
record when the collected works of M\"obius were ultimately
published.) The debate stopped there, apparently with the conclusion
that~$0^0$ should be undefined.

But no, no, ten thousand times no! Anybody who wants the binomial
theorem
$$(x+y)^n=\sum_{k=0}^n\,{n\choose k}\,x^ky^{n-k}\eqno(1.18)$$
to hold for at least one nonnegative integer $n$ {\it must\/} believe
that $0^0=1$, for we can plug in $x=0$ and $y=1$ to get~1 on the left
and~$0^0$ on the right. 

The number of mappings from the empty set to the empty set is~$0^0$.
It {\it has\/} to be~1.

On the other hand, Cauchy had good reason to
consider~$0^0$ as an undefined {\it limiting form}, in the sense that
the limiting value of $f(x)^{g(x)}$ is not known {\it a priori\/} when
$f(x)$ and $g(x)$  approach~0 independently. In this much stronger
sense, the value of~$0^0$ is less defined than, say, the value of
$0+0$. Both Cauchy and Libri were right, but Libri and his defenders
did not understand why truth was on their side.

Well, it's instructive to study mathematical history and to observe
how tastes change as progress is made. But let's come closer to the
present,
to see how Iverson's convention might be useful nowadays. Today's
mathematical literature is, in fact, filled with instances where
analogs of 
Iversonian brackets are being used---but the concepts must
be expressed in a roundabout way, because his convention is not yet
established. Here are two examples that I~happened to notice just
before writing this paper:

(1) Hardy and Wright, in the course of proving the Staudt-Clausen
theorem about the denominators of Bernoulli numbers~[20, \S~7.9],
consider the sum
$$\sum_{p-1\;{\rm divides}\;k}\;{1\over p}$$
where $p$ runs through primes. They define $\epsilon_k(p)$ to be 1 if
$p-1$ divides~$k$, otherwise $\epsilon_k(p)=0$; then the sum becomes
$$\sum_p{\epsilon_k(p)\over p}\,.$$
They proceed to show that $\sum_{m=1}^{p-1}m^k\equiv
-\epsilon_k(p)\pmod{p}$ whenever $p$ is prime, and the theorem follows
with a bit more manipulation.

(2) Mark Kac, introducing the relation of ergodic theory to continued
fractions~[24, \S~5.4], says: ``Let now $P_0\in\Omega$ and $g(P)$ the
characteristic function of the measurable set~$A$; i.e.,
$$g(P)=\cases{1,&$p\in A$,\cr
\noalign{\smallskip}
0,&$p\mathrel{\overline{\in}}A$.\cr}$$
It is now clear that $t(\tau,P_0,A)$ is given by the formula
$$t(\tau,P_0,A)=\int_0^{\tau}g\bigl(T_t(P_0)\bigr)\,dt\,,$$
and \dots~''.

I hope it is now clear why my students and I would find it quite
natural to say directly that
$$t(\tau,P_0,A)=\int_0^{\tau}[T_t(P_0)\in A]\,dt\,.$$
Also, in the context of Hardy and Wright, we would evaluate
$\biggl(\sum_{m=1}^{p-1}m^k\biggr)\bmod p$ 
and discover that it is $(p-1)\,[p-1\;{\rm divides}\;k]$. 

If you are a typical hard-working, conscientious mathematician,
interested in clear exposition and sound reasoning---and I~like to
include myself as a member of that set---then your experiences with
Iverson's convention may well go through several stages, just as mine
did. First, I~learned about the idea, and it certainly seemed
straightforward enough. Second, I~decided to use it informally while
solving problems. At this stage it seemed too easy to write just
`$[k\geq 0]$'; my natural tendency was to write something like
`$\delta(k\geq 0)$', giving an implicit bow to Kronecker, or
`$\tau(k\geq 0)$' where $\tau$~stands for truth. 
Adriano Garsia, similarly, decided to write `$\chi(k\geq 0)$', knowing
that $\chi$ often denotes a characteristic function; he has
used $\chi$ notation effectively in dozens of papers, beginning with
[10], and quite a few other mathematicians have begun to follow
his lead. (Garsia was one of my professors in graduate school, and
I~recently showed him the first draft of this note. He replied, ``My
definition from the very start was
$$\chi({\cal A})=\cases{1&if ${\cal A}$ is true\cr 
0&if ${\cal A}$ is false\cr}$$
where ${\cal A}$ is any statement whatever. But just like you, I~got
it by generalizing from Iverson's APL. \dots~I~don't have to tell you
the magic that the use of the $\chi$~notation can do.'')

If you go through the stages I did, however, you'll soon tire of
writing~$\delta$, $\tau$, or~$\chi$, when you recognize that the
notation is quite unambiguous without an additional symbol. 
Then you will have arrived
at the philosophical position adopted by Iverson when he wrote~[21].
And I~had also reached that stage when I~completed the first edition
of~[15]; I~adopted Iverson's original suggestion to enclose logical
statements in ordinary parentheses, not square brackets.

Unfortunately,  not all was  well with that first edition.
 Students found cases where I~had parenthesized a complicated
logical statement for clarity, for example when I~wrote something of
the form `$\alpha$ and $(\beta$ or~$\gamma)$'; they pointed out that
the simple act of putting parentheses around `$\beta$ or~$\gamma$'
automatically caused it to be evaluated as either~0 or~1, according to
a strict interpretation of Iverson's rule as I~had extended~it.

Worse yet, as I began to read the first edition of [15] with fresh
eyes, I~found that the formulas involved too many parentheses. It was
hard for me to perceive the structure of complex expressions that
involved Iversonian statements; the statements had been clear to me
when I~wrote them down, but they looked confusing when I~came back to
them several months later. A~computer could readily parse each
expression, but good notation must  be engineered for human beings.

Therefore in the second and subsequent printings of~[15], my
co-authors and~I now use square brackets instead of
parentheses, whenever we wish to transform logical statements into the
values~0 or~1. This resolves both problems, and we now believe that
the notation has proved itself well enough to be thrust upon the
world. Square brackets are used also for other purposes, but not in a
conflicting way, and not so often that the multiple uses become confusing.

One small glitch remains: We want to be able to write things like
$$\sum_p\,[p\;{\rm prime}]\,[p\leq x]/p\eqno(1.19)$$
to denote the sum of all reciprocals of primes $\leq x$. But this
summand unfortunately reduces to $0/0$ when $p=0$. In general, when an
Iverson-bracketed statement is false, we want it to evaluate into a
``very strong~0,'' namely a zero so strong that it annihilates
anything it is multiplied~by---even if that other factor is undefined.

Similarly, in formulas like (1.2) it is convenient to regard
${n\choose k}$ as strongly zero when $k$ is negative, so that, for
example, ${n\choose -10}z^{-10}=0$ when $z=0$. 

The strong-zero convention is enough to handle 99\%\ of the
difficult situations,
but we may also be using $1-[P(k)]$ to stand for the
quantity
[not~$P(k)$]; then we want $[P(k)]$ to give a ``strong~1.''  And
paradoxes can still arise, whenever irresistible forces meet immovable
objects. (What happens if a strong zero appears in the denominator?
And so~on.)

In spite of these potential problems in extreme cases, Iverson's
convention works beautifully in the vast majority of applications. It
is, in fact, far less dangerous than most of the other notations of
mathematics, whose dark corners we have learned to avoid long ago.
The safe use of Iverson's simple and convenient idea is quite easy to
learn. 

\bigskip\noindent
{\bf 2. Stirling numbers.}\enspace
The second plea I wish to make for perspicuous notation concerns the
famous coefficients introduced by James Stirling at the beginning of
his {\sl Methodus Differentialis\/} in 1730 [52]. The lack of a widely
accepted way to refer to these numbers has become almost scandalous. For
example, Goldberg, Newman, and Haynsworth begin their chapter on
Combinatorial Analysis in the NBS Handbook~[1] by remarking that
notations for Stirling numbers ``have never been standardized \dots
We feel that a capital~S is natural for Stirling numbers of the first
kind; it is infrequently used for other notation in this context. But
once it is used we have difficulty finding a suitable symbol for
Stirling numbers of the second kind. The numbers are sufficiently
important to warrant a special and easily recognizable symbol, and yet
that symbol must be easy to write. We have settled on a script
capital~{\eusm S} without any certainty that we have settled this
question permanently.''

The present predicament came about because Stirling numbers are indeed
important enough to have arisen in a wide variety of applications, yet
they are not quite important enough to have deserved a prominent place
in the most influential textbooks of mathematics. Therefore they have
been rediscovered many times, and each author has chosen a notation
that was optimized for one particular application.

The great utility of Stirling numbers has become clearer and
clearer with time, and mathematicians have now reached a stage where we can
intelligently choose a notation that will serve us well in the whole
range of applications.

I came into the picture rather late, having never heard of Stirling
numbers until after receiving my Ph.D. in mathematics. But I~soon
encountered them as I~was beginning to analyze the performance of
algorithms and to write the manuscript for my books on {\sl The Art of
Computer Programming}. I~quickly realized the truth of Imanuel Marx's
comment that ``these numbers have similarities with the binomial
coefficients~${n\choose k}$; indeed, formulas similar to those known
for the binomial coefficients are easily established''~[35]. In
order to emphasize those similarities and to facilitate pattern
recognition when manipulating formulas, Marx recommended using bracket
symbols~${n\brack k}$ for Stirling numbers of the first kind and brace
symbols~${n\brace k}$ for Stirling numbers of the second kind.
A~similar proposal was being made at about the same time in Italy by
Antonio Salmeri~[46].

I was strongly motivated by Charles Jordan's book, {\sl Calculus of
Finite Differences\/}~[23], which introduced me to the important
analogies between sums of factorial powers and integrals of ordinary
powers. 
But I~kept getting mixed up when I~tried to use Stirling numbers as he
defined them, because half of his ``first kind'' numbers were negative
and the other half were positive. I~had similar problems with Marx's
suggestions in~[35]; he made all Stirling numbers of the first kind
positive, but then he attached a minus sign to half the numbers of the
{\it second\/} kind. I~decided that I'd never be able to keep my head
above water unless I~worked with Stirling numbers that were entirely
signless.  

And I soon learned that the signless Stirling numbers have
important  combinatorial significance. So I~decided to try a
definition that combined the best qualities of the other notations
I'd~seen; I~defined the quantities ${n\brack k}$ and ${n\brace k}$
as follows:
$$\vcenter{\halign{\hfil$#\;$\hfil&#\hfil\cr
{n\brack k}=&the number of permutations of $n$ objects having $k$
cycles;\cr
\noalign{\smallskip}
{n\brace k}=&the number of partitions of $n$ objects into $k$ nonempty
subsets.\cr}}$$
For example, ${4\brack 2}=11$, because there are eleven different ways
to arrange four elements into two cycles:
$$\vcenter{\halign{$#$\hfil\qquad%
&$#$\hfil\qquad%
&$#$\hfil\qquad%
&$#$\hfil\cr
[1,2,3]\,[4]&[1,2,4]\,[3]&[1,3,4]\,[2]&[2,3,4]\,[1]\cr
[1,3,2]\,[4]&[1,4,2]\,[3]&[1,4,3]\,[2]&[2,4,3]\,[1]\cr
[1,2]\,[3,4]&[1,3]\,[2,4]&[1,4]\,[2,3].\cr}}$$
And ${4\brace 2}=7$, because the partitions of $\{1,2,3,4\}$ into two
subsets are
$$\vcenter{\halign{$#$\hfil\qquad%
&$#$\hfil\qquad%
&$#$\hfil\qquad%
&$#$\hfil\cr
\{1,2,3\}\{4\}&\{1,2,4\}\{3\}&\{1,3,4\}\{2\}&\{2,3,4\}\{1\}\cr
\{1,2\}\{3,4\}&\{1,3\}\{2,4\}&\{1,4\}\{2,3\}.\cr}}$$
Notice that this notation is mnemonic: The meaning of ${n\brace k}$
is easily remembered, because braces $\{\;\}$ are commonly used to
denote sets and subsets. We could also adopt the convention of
writing cycles in brackets, as in my examples above, where
$[1,2,3]=[2,3,1]=[3,1,2]$ is a typical three-cycle; that would make
the notation ${n\brack k}$ equally mnemonic. But I~don't insist on
this.

I have never decided how to pronounce `${n\brack k}$' and `${n\brace
k}$' when I'm reading formulas aloud in class. Many people have begun
to verbalize `${n\choose k}$' as ``$n$~choose~$k$''; hence I've been
saying ``$n$~cycle~$k$'' for ${n\brack k}$ and ``$n$~subset~$k$'' for
${n\brace k}$. But I~have also caught myself calling them
``$n$~bracket~$k$'' and ``$n$~brace~$k$.''

One of the advantages of these notational conventions is that binomial
coefficients and Stirling numbers can be defined by very simple
recurrence relations having a nice pattern:
$$\eqalignno{{n+1\choose k}&={n\choose k}+{n\choose k-1}\,;&(2.1)\cr
\noalign{\smallskip}
{n+1\brack k}&=n\,{n\brack k}+{n\brack k-1}\,;&(2.2)\cr
\noalign{\smallskip}
{n+1\brace k}&=k\,{n\brace k}+{n\brace k-1}\,.&(2.3)\cr}$$
Moreover---and this is extremely important---these identities hold for
all integers~$n$ and~$k$, whether positive, negative, or zero.
Therefore we can apply them in the midst of any formula (for example,
to ``absorb'' an~$n$ or a~$k$ that appears in the context
$n\,{n\brack k}$ or $k\,{n\brace k}$), without worrying about
exceptional circumstances of any kind.
 
I introduced these notations in the first  edition of my first
book~[25], and by now my students and~I have accumulated some
25~years of experience with them; the conventions have served us well.
However, such brackets and braces have still not become widely enough
adopted that they could be considered ``standard.'' For example,
Stanley's magnificent book on {\sl Enumerative Combinatorics\/}
[51] uses $c(n,k)$ for ${n\brack k}$ and $S(n,k)$ for~${n\brace
k}$. His notation conveys combinatorial significance, but it fails to
suggest the analogies to binomial coefficients that prove helpful in
manipulations. Such analogies were evidently not important enough in
his mind to warrant an extravagant two-line notation---although he
does use $\bigl(\!\bigl({n\atop k}\bigr)\!\bigr)$ to denote 
${n+k-1\choose k}=(-1)^k{-n\choose k}$, the number of combinations
with repetitions permitted. $\bigl($In a sense, Stanley's
$\bigl(\!\bigl({n\atop k}\bigr)\!\bigr)$ 
is a signless version of the numbers ${-n\choose k}$.$\bigr)$

When I wrote {\sl Concrete Mathematics\/} in 1988, I~explored Stirling
numbers more carefully than I~had ever done before, and I~learned two
things that really clinch the argument for~${n\brack k}$ and~${n\brace
k}$  as the best possible Stirling number notations. 
Ron Graham sent me a preview copy of a memorandum by B.~F. Logan
[34], which presented a number of interesting connections between
Stirling numbers and other mathematical quantities. One of the first
things that caught my attention was Logan's Table~1, a~two-dimensional
array that contained the numbers ${n\brack k}$ and ${n\brace k}$
simultaneously---implying that there really is only one ``kind'' of
Stirling number. Indeed, when I~translated Logan's results into my own
favorite notation, I~was astonished to find that his arrangement of
numbers was equivalent to a beautiful and easily remembered law of
duality, 
$${n\brace k}={-k\brack -n}\,.\eqno(2.4)$$
Once I had this clue, it was easy to check that the recurrence
relations (2.2) and (2.3) are equivalent to each other. And the
boundary conditions
$${0\brack k}={0\brace k}=[k=0]\qquad{\rm and}\qquad
{n\brack 0}={n\brace 0}=[n=0]\eqno(2.5)$$
yield unique solutions to (2.2) and (2.3) for all integers~$k$
and~$n$, when we run the recurrences forward and backward; the
``negative'' region for Stirling numbers of one kind turns out to
contain precisely the numbers of the other kind. For example, the following
subset of Logan's
 table gives the values of ${n\brack k}$ when $\vert n\vert$
and $\vert k\vert$ are at most~4:
$$\vcenter{\halign{%
$#\hfil$\qquad%
&$\hfil#\hfil$\quad%
&$\hfil#\hfil$\quad%
&$\hfil#\hfil$\quad%
&$\hfil#\hfil$\quad%
&$\hfil#\hfil$\quad%
&$\hfil#\hfil$\quad%
&$\hfil#\hfil$\quad%
&$\hfil#\hfil$\quad%
&$\hfil#\hfil$\cr
&k=-4&k=-3&k=-2&k=-1&k=0&k=1&k=2&k=3&k=4\cr
\noalign{\smallskip}
n=-4&1&0&0&0&0&0&\phantom{1}0&0&0\cr
n=-3&6&1&0&0&0&0&\phantom{1}0&0&0\cr
n=-2&7&3&1&0&0&0&\phantom{1}0&0&0\cr
n=-1&1&1&1&1&0&0&\phantom{1}0&0&0\cr
n=0&0&0&0&0&1&0&\phantom{1}0&0&0\cr
n=1&0&0&0&0&0&1&\phantom{1}0&0&0\cr
n=2&0&0&0&0&0&1&\phantom{1}1&0&0\cr
n=3&0&0&0&0&0&2&\phantom{1}3&1&0\cr
n=4&0&0&0&0&0&6&11&6&1\cr}}$$
The reflection of this matrix about a $45^{\circ}$ diagonal gives the
values of ${n\brace k}={-k\brack -n}$.

Naturally I wondered how I could have been working with Stirling
numbers for so many years without having been aware of such a basic
fact. Surely it must have been known before? After several hours of
searching in the library, I~learned that identity (2.4) had indeed been
known, but largely forgotten 
by succeeding generations of mathematicians,
 primarily because previous notations
for Stirling numbers made it impossible to state the identity in such
a memorable form. These investigations also turned up several things
about the history of Stirling numbers that I~had not previously
realized. 

During the nineteenth century, Stirling's connection with these
numbers had been almost entirely
forgotten. The numbers themselves were studied, in
the role of ``sums of products of combinations of the numbers
$\{1,2,\ldots,n\}$ taken~$k$ at a time.'' Let $C_k(n)$ and
$\Gamma_k(n)$ denote those sums, when the combinations are
respectively without or with repetitions; thus, for example,
$$\eqalign{C_4(4)&=1\cdot 2\cdot 3+1\cdot 2\cdot 4+1\cdot 3\cdot
4+2\cdot 3\cdot 4=50\,;\cr
\Gamma_3(3)&=1\cdot 1\cdot 1+1\cdot 1\cdot 2+1\cdot 1\cdot 3+1\cdot
2\cdot 2+1\cdot 2\cdot 3\cr
&\qquad\null+1\cdot 3\cdot 3+2\cdot 2\cdot 2+2\cdot 2\cdot 3+2\cdot
3\cdot 3+3\cdot 3\cdot 3=90\,.\cr}$$
It turns out that
$$C_k(n)={n+1\brack n+1-k}\qquad{\rm and}\qquad \Gamma_k(n)={n+k\brace
n}\,.\eqno(2.6)$$
Christian Kramp~[28]
proved near the end of the eighteenth century that
$$\eqalignno{C_k(n)&=\sum{n+1\choose k+l}\,{(k+l)!\over
j_1!\,2^{j_1}\,j_2!\,3^{j_2}\,j_3!\,4^{j_3}\,\ldots}\;,&(2.7)\cr
\noalign{\smallskip}
\Gamma_k(n)&=\sum{n+k\choose k+l}\,{(k+l)!\over
j_1!\,2!^{j_1}\,j_2!\,3!^{j_2}\,j_3!\,4!^{j_3}\,\ldots}\;,&(2.8)\cr}$$
where the sums are over all sequences of nonnegative integers $\langle
j_1,j_2,j_3,\ldots\,\rangle$ such that we have $j_1+2j_2+3j_3+\cdots =k$
(i.e., over all partitions of~$k$), and where
$l=j_1+j_2+j_3+\cdots\,\,$. For example, 
$$C_2(n)={n+1\choose 4}\,{1\over 8}+{n+1\choose 3}\,{1\over
3}\,;\qquad
\Gamma_2(n)={n+2\choose 4}\,{1\over 8}+{n+2\choose 3}\,{1\over
6}\,.$$
Notice that $C_k(n)$ and $\Gamma_k(n)$ are polynomials in~$n$, of
degree~$2k$.
The duality law (2.4) and the notational transformations of (2.6) are
equivalent to the amazing polynomial identity
$$C_k(n-1)=\Gamma_k(-n)\,;\eqno(2.9)$$
but hardly anybody was aware of this surpising fact,
otherwise we would almost certainly find it mentioned explicitly in
the comprehensive surveys compiled in the 1890s [19,~38].

On the other hand, a rereading of Stirling's original treatment~[52]
makes it clear that Stirling himself would not have found the duality law
(2.4) at all surprising. From the very beginning, he thought of
the numbers as two triangles hooked together in tandem. Indeed, his
entire motivation for studying them was the general identity
$$z^n=\sum_k{n\brace k}\,z^{\underline{k}}\,,\eqno(2.10)$$
which expresses ordinary powers in terms of falling factorial powers. 
When $n$ is positive, the nonzero terms in this sum occur for positive
values of $k\leq n$; but when $n$ is negative, the nonzero terms occur
for negative $k\leq n$. Stirling presented his tables by displaying
${n\brace k}$ with~$k$ as the row index and ${n\brack k}$ with~$k$ as
the column index; thus, he visualized a tandem arrangement exactly as
in the matrix of numbers above, with each column containing a sequence of
coefficients for (2.10).

I need to digress a bit about factorial powers. If $n$ is a positive
integer and $z$ is a complex number, I~like to write
$$z^{\underline{n}}=z(z-1)\,\ldots\,(z-n+1)\,,\eqno(2.11)$$
which I call ``$z$~to the~$n$ falling,'' and
$$z^{\overline{n}}=z(z+1)\,\ldots\,(z+n-1)\,,\eqno(2.12)$$
which is ``$z$ to the~$n$ rising.'' 
More generally, if $\alpha$ is any complex number, factorial powers
are defined by 
$$z^{\underline{\alpha}}=z!/(z-\alpha)!\qquad{\rm and}\qquad
z^{\overline{\alpha}}=\Gamma(z+\alpha)/\Gamma(z)\,,\eqno(2.13)$$
unless these formulas reduce to $\infty/\infty$ 
(when limiting values are used).
My use of underlined and overlined exponents is still controversial,
but I~cannot resist mentioning a curious fact: Many people (e.g.,
specialists in hypergeometric series) have become accustomed to the
notation $(z)_n$ for rising factorial powers, while many other people
(e.g., statisticians) use the same notation for {\it falling\/}
powers. The curious fact is that 
 this notation is called ``Pochhammer's symbol,'' but
Pochhammer himself~[43] used $(z)_n$ to stand for the binomial
coefficient~${z\choose n}$. I~prefer the underline/overline notation
because it is unambiguous and mnemonic, especially when I'm doing work
that involves factorial powers of both kinds. (Moreover, I~know
that~$z^{\underline{n}}$ and~$z^{\overline{n}}$ are easy to typeset,
using macros available in the file {\tt gkpmac.tex} in the standard
UNIX distribution of \TeX.)

In the special case $n=3$, Stirling's formula (2.10) gives
$$z^3={3\brace 3}\,z^{\underline{3}}+{3\brace 2}\,
z^{\underline{2}}+{3\brace 1}\,z^{\underline{1}}=
z(z-1)(z-2)+3z(z-1)+z\,.$$
And in the special case $n=-1$, it reduces to the infinite sum
$$\eqalignno{{1\over z}&=\sum_k{-1\brace k}\,z^{\underline{k}}\cr
\noalign{\smallskip}
&=\sum_k{\,k\,\brack 1}\,z^{\underline{-k}}\cr
\noalign{\smallskip}
&={0!\over z+1}+{1!\over (z+1)(z+2)}+{2!\over
(z+1)(z+2)(z+3)}+\cdots\,,&(2.14)\cr}$$
because
$${n\brack 1}=(n-1)!\,[n>0]\,.\eqno(2.15)$$
$\bigl($Stirling did not discuss convergence; he was, after all, writing in
1730. 
%It can be shown [NN2] that the infinite series (2.14) converges
%if and only if ${\rm Re}(z)>0$, and the same condition is necessary
%and sufficient for (2.10) when $n$ is a negative integer.)
We have the partial sum
$${1\over z}=\sum_{k=1}^n\,{(k-1)!\over (z+1)\,\ldots\,(z+k)}+
{n!\over z(z+1)\,\ldots\,(z+n)}\,;$$
this is a special case of the general identity
$${1\over z}=\sum_{k=1}^n\,{z_1\,\ldots\,z_{k-1}\over
(z+z_1)\,\ldots\,
(z+z_k)}+{z_1\,\ldots\,z_n\over
z(z+z_1)\,\ldots\,(z+z_n)}\eqno(2.16)$$
discovered by Fran{\c c}ois Nicole~[39] a few years before Stirling's
treatise appeared. Therefore the infinite series (2.14) converges if
and only if ${\rm Re}(z)>0$. By induction on~$n$, the same condition
is necessary and sufficient for (2.10) when~$n$ is any negative
integer. See~[41, \S$\,$30] for further discussion of (2.10).$\bigr)$

We noted above that the numbers ${m\brack k}$ can be regarded as sums
of products of combinations. The first identity in (2.6) is equivalent
to the formula
$$z^{\overline{n}}=\sum_k{n\brack k}\,z^k\,,\eqno(2.17)$$
when $n$
 is a nonnegative integer, if we expand the product~$z^{\overline{n}}$
and sum the coefficients of each power of~$z$.  Similarly, we have
$$z^{\underline{n}}=\sum_k{n\brack k}\,(-1)^{n-k}\,z^k\,.\eqno(2.18)$$
These equations are valid also when $n$ is a negative integer; in that case
both infinite series converge for $\vert z\vert>\vert n\vert$. Notice that
(2.10) and (2.18) tell us how to convert back and forth between
ordinary powers and factorial powers.

Let's turn now to the nineteenth century.
 Kramp~[29] decided to explore a slightly generalized type
of factorial power, for which he used the notations
$$\eqalignno{a^{n\mid
r}&=a(a+r)\,\ldots\,\bigl(a+(n-1)\,r\bigr)&(2.19)\cr
\noalign{\smallskip}
a^{-n\mid r}&=1/(a-r)(a-2r)\,\ldots\,(a-nr)&(2.20)\cr}$$
when $n$ is a positive integer. Then he considered the expansion
$$a^{n\mid r}=a^n+n\,\hbox{\eufr k}\,1.\,a^{n-1}r+n\,\hbox{\eufr k}\,
2.\,
a^{n-2}r^2+\cdots \;,\eqno(2.21)$$
where the coefficients $n\,\hbox{\eufr k}\,m$ are independent of~$a$
and~$r$ [29, \S\S 539--540]; thus, $n\,\hbox{\eufr k}\,m$ was his
notation for ${n\brack n-m}$. He obtained~[29, \S~557] a~series of
formulas equivalent to
$$m\,{n\brack n-m}=\sum_{k=0}^{m-1}{n-k\choose m+1-k}\,{n\brack
n-k}\,,\eqno(2.22)$$
thereby giving a new proof that ${n\brack n-m}$ is a polynomial in~$n$
of degree~$2m$. This proof, independent of his earlier formulas~(2.7)
and~(2.8), works for both positive and negative values of~$n$.

Kramp implicitly understood the duality principle (2.4), in the sense
that he regarded the coefficients~${\,n\,\brack k}$ and~${n\brace k}$
as the positive and negative portions of a doubly infinite array  of
numbers. In fact, he assumed that equation (2.21) would hold for
arbitrary real values of~$n$. He differentiated $a^{x\mid r}$ with
respect to~$x$ and gave formal derivations of several interesting
series. However, his expansion (2.21) is equivalent~to
$$z^{\overline{n}}=\sum_k\,{n\brack n-k}\,z^{n-k}\eqno(2.23)$$ 
\noindent
$\bigl($a slight variation of (2.17)$\bigr)$, and this series is not
always convergent for noninteger~$n$.
We can show, for example, that
$$\Biggl\vert{1/2\brack 1/2-k}\Biggr\vert 
>k!/7^k\qquad\hbox{for infinitely many
$k$}\,;\eqno(2.24)$$
hence (2.23) diverges for all $z$ when $n=1/2$. Kramp lived before the
days when convergence of infinite series was understood. (See~[29, \S~574],
where he says that the divergent series $\sum_{k>0}B_ky^k\!/k$ is
``tr\`es convergente pour peu que~$y$ soit une petite fraction''!)

Several other nineteenth-century authors developed the theory of
factorial powers, notably Andreas von Ettingshausen~[6], Ludwig
Schl\"afli [41, 48], and Oskar Schl\"omilch [49], who used the
respective notations
$$\mathop{\rm F\null}^n\!{}_m\,,\qquad
\mathop{A\null}^n\!{}_m\,,\qquad
{\rm and}\qquad
\mathop{C\null}^n\!{}_m$$
for the coefficients ${n\brack n-m}$. All of these authors considered
both positive and negative integers~$n$. Thus, for example,
Ettingshausen's notation for a Stirling number such as ${n+m\brace n}$
$={-n\brack -n-m}$ was 
$$\mathop{\rm{F}\rlap{$_m$}}^{-n}$$ (see~[6, \S~151]).

Incidentally, these works of Kramp and Ettingshausen proved to be
important in the history of mathematical notations. Kramp's book
introduced the notation~$n!$ for factorials~[29, pages~V and 219], and
Ettingshausen's book introduced the notation~${n\choose k}$ for
binomial coefficients~[6, page~30]. Ettingshausen wrote his book
shortly after Fourier~[8] had invented $\sum$-notation for sums;
Ettings\-hausen tried a German variation, writing
$\hbox{\eufr S}^k_{a,b}$ for what has evolved into $\sum_{k=a}^b$.
%$\hbox{\eufr S}\limits^k_{a,b}$ for what has evolved into $\sum_{k=a}^b$.
He also wrote $(a,r)^n$ for Kramp's $a^{n\mid r}$; thus, for example,
Ettingshausen~[6, \S~153 and \S~156] gave the equations
$$(a,d)^n=\mathop{\hbox{\eufr S}}^w_0\,\mathop{\rm F\null}^n\!{}\!_w\,
a^{n-w}\,d^w\qquad{\rm and}\qquad
a^n=\mathop{\hbox{\eufr S}}^r_0\,(-1)^r\,
\mathop{\rm F\rlap{$_r$}}^{-n+r}\,(a,d)^{n-r}\,d^r$$ 
as equivalents of Kramp's (2.21) and Stirling's (2.10). He presented
Kramp's (2.22) in the form
$$v\,
\mathop{\rm F\null}^n\!{}\!_v
=\mathop{\hbox{\eufr S}}^w_{0,v-1}\,{n-w\choose v+1-w}\,
\mathop{\rm F\null}^n\!{}\!_w\,,$$
and remarked~[6, \S~154] that this holds for both negative and
positive~$n$. Ettingshausen had related the F~coefficients
to sums of products of combinations
with and without repetition; thus he implicitly confirmed (2.9).

The first person to attach Stirling's name to the numbers we now
call Stirling numbers was Niels Nielsen in 1904~[40]; he said that
this new nomenclature had been suggested to him by T.~N. Thiele.
(The numbers may have been studied before Stirling's time; for
example, I~once found the values of~${n\brack k}$ for $1\leq n\leq 7$
in some unpublished manuscripts of Thomas Harriot, dating from about
1600, in the British Museum~[26, page~241]. But Stirling almost
surely deserves the credit for being first to deduce nontrivial facts
about ${n\brack k}$ and~${n\brace k}$.)

Nielsen wrote~$C_n^k$ for ${n\brack n-k}$, which he called a
``Stirling number of rank~$n$''; and he wrote $\hbox{\eufr C}_n^k$
for ${n+k-1\brace n-1}$, which he called a ``Stirling number of
rank~$-n$.'' (He should really have defined its rank to be $1-n$).
In equation (41) of his paper, Nielsen obtained a rigorous proof
of the duality law (2.4); but he had to state it in a peculiar way,
because he had defined $C_n^k$ and $\hbox{\eufr C}_n^k$ only for
nonnegative~$n$ and~$k$. Thus, he could not write $C_n^k=\hbox{\eufr
C}_{1-n}^k$; he had to say instead that $f_k(n)=g_k(1-n)$, where
$f_k(n)$ and $g_k(n)$ were the polynomials defined by~$C_n^k$
and~$\hbox{\eufr C}_n^k$. Tweedie [54] expressed (2.4) with similar
circumlocutions. 

When Jordan took up Stirling numbers [22], he wrote $S_n^k$ for
$(-1)^{n-k}{n\brack k}$ and $\hbox{\eufr S}_n^k$ for~${n\brace k}$. He
does not seem to have known the duality law (2.4), probably because he
had learned about Stirling numbers from Nielsen's book~[41], which
omitted some of the details in Nielsen's paper~[40].
And as far as I~know, the duality law largely disappeared from
mathematicians' collective consciousness during most of the twentieth
century; it seems to have been mentioned explicitly
only in a few scattered places: (1)~Hansraj Gupta, ``working in a
small township away from what was then the only University in the
Panjab'' 
[18, page~5], rediscovered Stirling numbers and Stirling duality by
himself, in the early 1930s. This became part of his Ph.D.
dissertation [17], and he included it in a book on number theory prepared
many years later [18, Chapter~5]. (2)~H.~W. Gould~[12] 
was probably the first  twentieth-century mathematician
to observe that we can use the
polynomials ${n\brack n-k}$ and ${n\brace n-k}$ to extend the domain
of Stirling numbers to negative values of~$n$. Gould's way of writing
(2.4) was $S_1(-n-1,k)=S_2(n,k)$; and shortly thereafter~[13], he
mentioned the equivalent formula
$$S_{-k}^{-n}=(-1)^{n-k}\hbox{\eufr S}_n^k\,,$$
in Jordan's notation. 
(3)~R.~V. Parker [42], like Gupta, displayed both of Stirling's
triangles in tandem, presenting them in a single table as Logan later
did. (4)~In 1976, Ira Gessel and Richard Stanley investigated some of
the deeper structure underlying the Stirling polynomials
$f_k(n)={n+k\brace n}$ and $g_k(n)={n\brack n-k}$. They noted in
particular [11, equation~(3)] that $f_k(-n)=g_k(n)$. This fact is
equivalent to the duality law (2.4).

Stanley had discovered a beautiful theorem in his Ph.D. thesis a few
years earlier [50, Propostion 13.2(i)], now called the reciprocity
theorem for order polynomials: If $P$ is any finite partially ordered
set, let $\Omega(P,n)$ be the number of order-preserving mappings
from~$P$ into the totally ordered set $\{1,2,\ldots,n\}$; and let
$\overline{\Omega}(P,n)$ be the number of such mappings that are
strictly order-preserving. Thus, if $x\prec y$ in~$P$, the
mappings~$f$ enumerated by $\Omega(P,n)$ must satisfy $f(x)\leq f(y)$,
and the mappings~$g$ enumerated by $\overline{\Omega}(P,n)$ must
satisfy $g(x)<g(y)$. Stanley's theorem states that, in general, we
have $f(-n)=(-1)^pg(n)$, where $p$ is the number of elements of~$P$.
For example, if $P$ consists of $p$~isolated points with no order
constraints whatever, we have
$\Omega(P,n)=\overline{\Omega}(P,n)=n^p$. And if the points of~$P$ are
themselves totally ordered, then $\Omega(P,n)$ is ${n+p-1\choose p}$,
the number of combinations of $n$ things $p$ at a time with
repetitions permitted, and $\overline{\Omega}(P,n)$ is ${n\choose p}$,
the combinations without repetition. In both cases we have
$\Omega(P,-n)=(-1)^p\,\overline{\Omega}(P,n)$. 

I showed Stanley the first draft of this note and asked him whether
the Stirling duality law (2.4) could be derived as a special case of
his general reciprocity law. Sure enough, he replied that Gessel had
noticed a simple way to do exactly that, shortly after the paper [11]
was written. Let $P_k$ be the partial order on $2k$~points typified by
$$\unitlength=15pt
\def\\(#1){\put(#1){\disk{.1}}}
P_4=
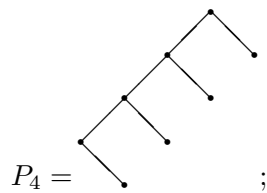
\beginpicture(4,4)(0,0)
\\(0,1)\\(1,0)\\(1,2)\\(2,1)\\(2,3)\\(3,2)\\(3,4)\\(4,3)
\put(0,1){\line(1,1)3}
\def\\(#1){\put(#1){\line(-1,1)1}}%
\\(1,0)\\(2,1)\\(3,2)\\(4,3)
\endpicture
\,;$$
then
$$\eqalign{\Omega(P_k,n)&=\sum_{1\leq
x_1,\ldots,x_k,y_1,\ldots,y_k\leq n}[x_1\leq \cdots \leq x_k][x_1\geq
y_1]\,\ldots\,[x_k\geq y_k]\cr
\noalign{\smallskip}
&=\sum_{1\leq x_1,\ldots,x_k\leq n}[x_1\leq \cdots \leq x_k]
\,x_1\,\ldots\,x_k\,,\cr}$$
and
$$\eqalign{\overline{\Omega}(P_k,n)&=\sum_{1\leq
x_1,\ldots,x_k,y_1,\ldots,y_k\leq n}[x_1< \cdots<  x_k][x_1 >
y_1]\,\ldots\,[x_k > y_k]\cr
\noalign{\smallskip}
&=\sum_{2\leq x_1,\ldots,x_k\leq n}[x_1 < \cdots <
x_k](x_1-1)\,\ldots\,(x_k-1) \cr
\noalign{\smallskip}
&=\sum_{1\leq x_1,\ldots,x_k\leq n-1}[x_1 < \cdots <
x_k]\,x_1\,\ldots\,x_k
\,. \cr}$$
Thus the sums are respectively $\Gamma_k(n)$ and $C_k(n-1)$; by (2.6)
we have $\Omega(P_k,n)={n+k\brace n}$ and
$\overline{\Omega}(P_k,n)={n\brack n-k}$, hence (2.4) is indeed an
instance of Stanley's theorem.

Now we are ready to discuss the second reason why I became convinced that
${n\brack k}$ is the right symbolism for these coefficients
after I~had translated Logan's memo [34] into that notation: We know
that ${n\brack n-k}$ is a polynomial in~$n$, when $k$
 is an integer; hence, as Kramp knew,  we can sensibly define the quantity
${\alpha\brack \alpha-k}$ for arbitrary complex~$\alpha$ and
integer~$k$, using that same polynomial. Then---and here comes the
punch line---Logan noticed that the fundamental 
equations (2.17) and
(2.18) generalize to {\it asymptotic formulas}, valid for arbitrary
exponents~$\alpha$: If $z\rightarrow\infty$ and if $m$~is any
nonnegative integer, we have
$$\eqalignno{z^{\overline{\alpha}}&=\sum_{k=0}^m{\alpha\brack
\alpha-k}\,z^{\alpha-k}+O(z^{\alpha-m-1})\,;&(2.25)\cr
\noalign{\smallskip}
z^{\underline{\alpha}}&=\sum_{k=0}^m{\alpha\brack
\alpha-k}\,(-1)^k\,z^{\alpha-k}+O(z^{\alpha-m-1})\,.&(2.26)\cr}$$
(See~[15, exercise 9.44]; 
equation (2.25) is a correct way to formulate Kramp's divergent series
(2.23). These equations are special cases of a
still more general result proved by 
Tricomi and Erd\'elyi [53,~9].) The easily remembered
expansions in
 (2.25) and (2.26) were quite a revelation to me. I~had
often spent time laboriously calculating approximations to ratios such
as $z^{\overline{1/2}}=\Gamma(z+1/2)/\Gamma(z)$, the hard way: I~took
logarithms, then used Stirling's approximation, and then took
exponentials. But equations (2.25) and (2.26) produce the answer directly.

Moreover Stirling's original identity (2.10) can be generalized in a
similar way: If $\alpha$ is any complex number, we have
$$z^{\alpha}=\sum_k\,{\alpha\brace\alpha-k}\,z^{\underline{\alpha-k}}\,,
\qquad {\rm Re}(z)>0\,.\eqno(2.27)$$
When I wrote the first draft of this note, I~knew only that the series
(2.27) was convergent, and that it was asymptotically correct as
$z\rightarrow\infty$; so I~conjectured that equality might hold. Soon
afterward, B.~F. Logan found the following proof (although he
naturally stated
it in his own notation): Suppose first that Re$(\alpha)<1$. Then we
have the well known identity
$$z^{\alpha-1}={1\over\Gamma(1-\alpha)}\,
\int_0^{\infty}e^{-zt}t^{-\alpha}\,dt\,,\qquad{\rm
Re}(z)>0\,,\eqno(2.28)$$
and we can substitute $e^{-t}=1-u$ to get
$$z^{\alpha-1}={1\over\Gamma(1-\alpha)}\,\int_0^1(1-u)^{z-1}u^{-\alpha}
\left(\,{1\over
u}\,\ln\,{1\over 1-u}\,\right)^{-\alpha}\,du\,.$$
Now it turns out that the powers of ${1\over u}\,\ln\,{1\over 1-u}$
generate the Stirling numbers
${\alpha\brace\alpha-k}={k-\alpha\brack-\alpha}$, in the sense that
$$\left(\,{1\over u}\,\ln\,{1\over
1-u}\,\right)^{-\alpha}=\sum_k\,{\alpha\brace\alpha-k}\,{u^k\over
(k-\alpha)\,\ldots\,(1-\alpha)}\,,\eqno(2.29)$$
a series that converges for $\vert u\vert<1\;\;\bigl($see [15,
equations (6.45), (6.53), (7.50)]$\bigr)$. Therefore
$$\eqalign{z^{\alpha}&=\sum_k\,{\alpha\brace\alpha-k}\,{z\over
\Gamma(k+1-\alpha)}\,\int_0^1(1-u)^{z-1}u^{k-\alpha}\,du\cr
\noalign{\smallskip}
&=\sum_k\,{\alpha\brace\alpha-k}\,{\Gamma(z+1)\over\Gamma(z+1+k-\alpha)}
=\sum_k\,{\alpha\brace\alpha-k}\,{z!\over(z+k-\alpha)!}\,,\cr}$$
and (2.27) is verified when Re$(\alpha)<1$. To complete the proof,
we need only show that (2.27) holds for $\alpha+1$ if it holds
for~$\alpha$; but this is easy, because
$$\eqalign{z^{\alpha+1}&=\sum_k\,{\alpha\brace\alpha-k}\,z\cdot
z^{\underline{\alpha-k}}\cr 
\noalign{\smallskip}
&=\sum_k\,{\alpha\brace\alpha-k}\,\left(z^{\underline{\alpha+1-k}}+(\alpha
-k)z^{\underline{\alpha-k}}\right)\cr
\noalign{\smallskip}
&=\sum_k\,{\alpha\brace\alpha-k}\,z^{\underline{\alpha+1-k}}+
\sum_k\,{\alpha\brace\alpha+1-k}\,(\alpha+1-k)z^{\underline{\alpha+1-k}}\cr
\noalign{\smallskip}
&=\sum_k\,{\alpha+1\brace\alpha+1-k}\,z^{\underline{\alpha+1-k}}\cr}$$
by the basic recurrence equation (2.3).

Notice that in all of the general identities (2.25)--(2.27), 
as in the original formulas (2.10), (2.17), and (2.18) that inspired
them, the lower index
within the braces or brackets is the same as the exponent of~$z$. This
makes the relations easy to remember, by analogy with the binomial
theorem
$$(1+z)^{\alpha}=\sum_k{\alpha\choose k}\,z^k\,,\qquad{\rm when}\;
\vert z\vert <1\,.\eqno(2.30)$$

Some readers will have been thinking, ``This all looks fairly
plausible, but unfortunately Knuth is overlooking a key point that
ruins the whole proposal: We {\it can't\/} use the notation ${n\brack
k}$ for Stirling numbers, because it has already been used for more
than a century as the standard notation for Gauss's generalized
binomial coefficients.''

Well, there is a down side to every good idea, but this objection is
not really severe.
For one thing, the standard notation for Gaussian binomial
coefficients involves a hidden parameter~$q$, and it's not unusual for
modern researchers to make transformations that change~$q$. Therefore
Gauss's notation  is incomplete, and Andrews (for example) has used
the notation ${n\brack k}_{q^2}$ for the Gaussian coefficient 
with~$q^2$ as the hidden parameter~[2, page~49].
Such examples suggest that it is appropriate to denote
 Gaussian binomials as ${n\choose
k}_q$, especially since they reduce to ordinary binomials when $q=1$.
This notation also generalizes nicely to such things as Fibonomial
coefficients ${n\choose k}_{\hbox{\eusm F}}\,$; see~[27]. We can 
then reserve the notation ${n\brack k}_q$ for a
$q$--generalization of~${n\brack k}$. (The reverse strategy was
unfortunately adopted in~[14].) Secondly, I~do not believe that any
existing mathematical works, including  books
like~[2] which use Gaussian coefficients extensively, 
 would become seriously cluttered if the Gaussian ${n\brack
k}$ were changed everywhere to ${n\choose k}_q$. 
Even so, such changes are not necessary;
there is obviously no
harm in beginning a mathematical paper or a book chapter or an entire
book with a statement to the effect that ``${n\brack k}$ will denote a
Gaussian binomial coefficient 
with parameter~$q$ in what follows.'' All notation can be
redefined for special purposes. Therefore Stirling number enthusiasts
are not encroaching on Gaussian territory when they write ${n\brack
k}$, if they also mumble something about Stirling in order to set the
context.

One further point is worth noting in conclusion:
As soon as the notations ${n\brack k}$ and/or ${n\brace k}$ are
adopted, there will no longer be a need to speak about Stirling
numbers ``of the first and second kind,'' except as a concession to
history. Nielsen wrote a superb book~[41], but he did the world a
disservice by originating the {\it Erster Art\/} and {\it Zweier Art\/}
terminology, because that terminology  has no mnemonic value and 
 is historically
inaccurate. Stirling introduced the numbers ${n\brace k}$ first and
brought in ${n \brack k}$ second. Indeed, practical applications have
always tended to involve the numbers ${n\brace k}$ much more often than their
${n\brack k}$ counterparts. It seems far better to speak of ${n\brace
k}$ as a Stirling subset number, and to call ${n\brack k}$ a Stirling
cycle number. Then the names are tied to intuitive, student-friendly
concepts, not to arbitrary and offputting concepts of the $k$th kind.

\bigskip\noindent
{\bf Acknowledgments.}\enspace I am extremely grateful for comments
received from John Ewing, Philippe Flajolet, Adriano Garsia, B.~F.
Logan, Andrew Odlyzko, Richard Stanley, and H.~S. Wilf, without
which these notes would have been substantially poorer.

\bigskip
\centerline{\bf References}

\bib
\phantom{1}[1]\enspace
Milton Abramowitz and Irene A. Stegun, editors, {\sl Handbook of
Mathematical Functions\/} (U.S. National Bureau of Standards, 1964).

\bib
\phantom{1}[2]\enspace
George E. Andrews, {\sl The Theory of Partitions}, Encyclopedia of
Mathematics and its Applications, volume~2 (Reading, Mass.:
Addison\kern-.1em--Wesley, 1976).

\bib
\phantom{1}[3]\enspace 
Anonymous and S{\thinspace}\dots , 
``Bemerkungen zu den Aufsatze \"uberschrieben,
`Beweis der Gleichung $0^0=1$ nach J.~F. Pfaff,' im zweiten Hefte
dieses Bandes, S.~134,''
{\sl Journal f\"ur die reine und angewandte Mathematik\/
\bf 12} (1834), 292--294.

\bib
\phantom{1}[4]\enspace
Charles Babbage, {\sl Passages from the Life of a Philosopher\/}
(London, 1864). Reprinted in {\sl Charles Babbage and his Calculating
Engines}, edited by Philip Morrison and Emily Morrison (New York:
Dover, 1961).

\bib
\phantom{1}[5]\enspace
Augustin-Louis Cauchy, {\sl Cours d'Analyse de l'Ecole Royale
Polytechnique\/} (1821). In his {\sl {\OE}uvres Compl\`etes},
series~2, volume~3.

\bib
\phantom{1}[6]\enspace
Andreas v.\ Ettingshausen, {\sl Die combinatorische Analysis\/}
(Vienna, 1826).

\bib
\phantom{1}[7]
Philippe Flajolet and Andrew Odlyzko, ``Singularity analysis of
generating functions,'' {\sl SIAM Journal on Discrete Mathematics\/
\bf 3} (1990), 216--240.

\bib
\phantom{1}[8]\enspace
J. Fourier, ``Refroidissement s\'eculaire du globe terrestre,'' {\sl
Bulletin des Sciences par la Soci\'et\'e philomathique de Paris},
series~3, {\bf 7} (1820), 58--70. Reprinted in {\sl{\OE}uvres de
Fourier}, volume~2, 271--288.

\bib
\phantom{1}[9]\enspace
C. L. Frenzen, ``Error bounds for asymptotic expansions of the
ratio of two gamma functions,'' {\sl SIAM Journal on Mathematical
Analysis\/ \bf 18} (1987), 890--896.

\bib
[10]\enspace
Adriano M. Garsia, ``On the `maj' and `inv' $q$-analogues of Eulerian
polynomials,'' {\sl Linear and Multilinear Algebra\/ \bf 8} (1979),
21--34. 

\bib
[11]\enspace
Ira Gessel and Richard P. Stanley, ``Stirling polynomials,'' {\sl
Journal of Combinatorial Theory\/ \bf A24} (1978), 24--33.

\bib
[12]\enspace
H. W. Gould, ``Stirling number representation problems,'' {\sl
Proceedings of the American Mathematical Society\/ \bf 11} (1960),
447--451. For subsequent work, see his review of [42] in {\sl
Mathematical Reviews\/ \bf 49} (1975), 885--886.

\bib
[13]\enspace
H. W. Gould, ``Note on a paper of Klamkin concerning Stirling numbers,
This {\sc Monthly} {\bf 68} (1961), 477--479.

\bib
[14]\enspace
H. W. Gould, ``The $q$--Stirling numbers  of first and second kinds,''
{\sl Duke Mathematical Journal\/ \bf 28} (1961), 281--289.

\bib
[15]\enspace
Ronald L. Graham, Donald E. Knuth, and Oren Patashnik, {\sl Concrete
Mathematics\/} (Reading, Mass.: Addison\kern-.1em--Wesley, 1989).

\bib
[16]\enspace
Daniel H. Greene and Donald E. Knuth, {\sl Mathematics for the
Analysis of Algorithms}, second edition (Boston: Birkh\"auser, 1981).
Third edition, 1990.

\bib
[17]\enspace
H. Gupta, {\sl Symmetric Functions in the Theory of Integral Numbers},
Lucknow University Studies {\bf 14} (Allahabad: Allahabad Law Journal
Press, 1940).

\bib
[18]\enspace
Hansraj Gupta, {\sl Selected Topics in Number Theory\/} (Tunbridge
Wells, England: Abacus Press, 1980).

\bib
[19]\enspace
Johann G. Hagen, {\sl Synopsis der H\"oheren Mathematik\/ \bf 1}
(Berlin, 1891).

\bib
[20]\enspace
G. H. Hardy and E. M. Wright, {\sl An Introduction to the Theory of
Numbers\/} (Oxford, Clarendon Press, 1938). Fifth edition, 1979.

\bib
[21]\enspace
Kenneth E. Iverson, {\sl A Programming Language\/} (New York: Wiley,
1962). 

\bib
[22]\enspace
Charles Jordan, ``On Stirling's Numbers,'' {\sl T{\^o}hoku
Mathematical Journal\/ \bf 37} (1933), 254--278.

\bib
[23]\enspace
Charles Jordan, {\sl Calculus of Finite Differences}
(Budapest, 1939). Third edition, 1965.

\bib
[24]\enspace
Mark Kac, {\sl Statistical Independence in Probability, Analysis and
Number Theory}, Carus Mathematical Monographs {\bf 12} (Mathematical
Association of America, 1959).

\bib
[25]\enspace
Donald E. Knuth,
{\sl Fundamental Algorithms\/} (Reading, Mass.:
Addison\kern.1em--Wesley, 1968).

\bib
[26]\enspace
Donald E. Knuth, review of {\sl History of Binary and Other Nondecimal
Numeration\/} by Anton Glaser, {\sl Historia  Mathematica\/ \bf 10}
(1983), 236--243.

\bib
[27]\enspace
Donald E. Knuth and Herbert S. Wilf, ``The power of a prime that
divides a generalized binomial coefficient,'' {\sl Journal f\"ur die
reine und angewandte Mathematik\/ \bf 396} (1989), 212--219.

\bib
[28]\enspace
Christian Kramp, ``Coefficient des allgemeinen Gliedes jeder
willk\"uhrlichen Potenz eines Infinitinomiums; Verhalten zwischen
Coefficienten der Gleichungen und Summen der Produkte und der Potenzen
ihrer Wurzeln; Transformation und Substitution der Reihen durch
einander,'' in {\sl Der polynomische Lehrsatz}, edited by Carl
Friedrich Hindenburg (Leipzig, 1796), 91--122.

\bib
[29]\enspace
C. Kramp, {\sl \'El\'emens d'arithm\'etique universelle\/}
(Cologne, 1808).

\bib
[30]\enspace
Leopold Kronecker, ``Ueber bilineare Formen,''
{\sl Journal f\"ur die reine und angewandte Mathematik\/
\bf 68} (1868), 273--285.

\bib
[31]\enspace
Leopold Kronecker, {\sl Vorlesungen \"Uber de Theorie der
Determinanten}, edited by Kurt Hensel, volume~1 (Leipzig: Teubner, 1903).

\bib
[32]\enspace
Guillaume Libri, ``Note sur les valeurs de la fonction~$0^{0^x}$,''
{\sl Journal f\"ur die reine und angewandte Mathematik\/
\bf 6} (1830), 67--72.

\bib
[33]\enspace
Guillaume Libri, ``M\'emoire sur les fonctions discontinues,''
{\sl Journal f\"ur die reine und angewandte Mathematik\/
\bf 10} (1833), 303--316.

\bib
[34]\enspace
B. F. Logan, ``Polynomials related to the Stirling numbers,'' 
AT\&T Bell Labs internal technical memorandum, August~10, 1987.

\bib
[35]\enspace
Imanuel Marx, ``Transformation of series by a variant of Stirling
numbers,'' This {\sc Monthly} {\bf 69} (1962), 530--532. His ${n\brack
k}$ is my ${n+1\brack k+1}$; his ${n \brace k}$ is my
$(-1)^{n-k}{n+1\brace k+1}$.

\bib
[36]\enspace
A. F. M\"obius, ``Beweis der Gleichung $0^0=1$, nach J. F. Pfaff,''
{\sl Journal f\"ur die reine und angewandte Mathematik\/
\bf 12} (1834), 134--136.

\bib
[37]\enspace
Douglas H. Moore, {\sl Heaviside Operational Calculus: An Elementary
Foundation\/} (New York: American Elsevier, 1971).

\bib
[38]\enspace
Eugen Netto, {\sl Lehrbuch der Combinatorik\/} (Leipzig, 1901). Second
edition, with additions by Thoralf Skolem and Viggo Brun, 1927.

\bib
[39]\enspace
Nicole, ``M\'ethode pour sommer une infinit\'e de Suites nouvelles,
dont on ne peut trouver les Sommes par les M\'ethodes connu\"es,''
{\sl M\'emoires de l'Academie Royale des Sciences\/} (Paris, 1727),
257--268. 

\bib
[40]\enspace
Niels Nielsen, ``Recherches sur les polynomes et les nombres de
Stirling,''
{\sl Annali di Matematica pura ed applicata}, series~3, {\bf 10}
(1904), 287--318.

\bib
[41]\enspace
Niels Nielsen, {\sl Handbuch der Theorie der Gammafunktion\/}
(Leipzig: Teubner, 1906).

\bib
[42]\enspace
R. V. Parker, ``The complete polynomial grid,'' {\sl Matematichki
Vesnik\/ \bf 10} ({\bf 25}) (1973), 181--203.

\bib
[43]\enspace
L. Pochhammer, ``Ueber hypergeometrische Functionen $n^{\rm ter}$
Ordnung,'' 
{\sl Journal f\"ur die reine und angewandte Mathematik\/
\bf 71} (1870), 316--352.

\bib
[44]\enspace
Hillel Poritsky, ``Heaviside's operational calculus---its applications
and foundations,'' This {\sc Monthly}, {\bf 43} (1936), 331--344.

\bib
[45] \enspace
S{\thinspace}\dots , ``Sur la valeur de $0^0$,'' 
{\sl Journal f\"ur die reine und angewandte Mathematik\/
\bf 11} (1834), 272--273.

\bib
[46]\enspace
Antonio Salmeri, ``Introduzione alla teoria dei coefficienti
fattoriali,'' {\sl Giornale di Matematiche di Battaglini\/ \bf 90}
(1962), 44--54. His ${n\brack k}$ is my ${n+1\brack n+1-k}$.

\bib
[47]\enspace
Schlaeffli, ``Sur les co\"efficients du d\'eveloppement du produit
$1.(1+x)(1+2x)\,\ldots\,\bigl(1+(n-1)\,x\bigr)$ suivant les puissances
ascendantes de~$x$,'' 
{\sl Journal f\"ur die reine und angewandte Mathematik\/
\bf 43} (1852), 1--22.

\bib
[48]\enspace
Schl\"affli, ``Erg\"anzung der Abhandlung \"uber die Entwickelung des
Products 
$1.(1+x)(1+2x)$\allowbreak
$(1+3x)\,\ldots\,\bigl(1+(n-1)\,x\bigr)=
{\displaystyle\mathop{\Piit\null}\limits^n(x)}$
%\mathop{\mit\Pi\null}\limits^n(x)$
in Band XLIII dieses Journals,''
{\sl Journal f\"ur die reine und angewandte Mathematik\/
\bf 67} (1867), 179--182.

\bib
[49]\enspace
O. Schl\"omilch, ``Recherches sur les coefficients des facult\'es
analytiques,''
{\sl Journal f\"ur die reine und angewandte Mathematik\/
\bf 44} (1852), 344--355.

\bib
[50]\enspace
Richard P. Stanley, {\sl Ordered Structures and Partitions}, Memoirs
of the American Mathematical Society {\bf 119} (1972).

\bib
[51]\enspace
Richard P. Stanley, {\sl Enumerative Combinatorics}, volume~1
(Belmont, Calif.: Wadsworth, 1986).

\bib
[52]\enspace
James Stirling, {\sl Methodus Differentialis\/} (London, 1930).
English translation, {\sl The Differential Method}, 1749.

\bib
[53]\enspace
F. G. Tricomi and A. Erd\'elyi, ``The asymptotic expansion of a ratio
of gamma functions,'' {\sl Pacific Journal of Mathematics\/ \bf 1}
(1951), 133--142.

\bib
[54]\enspace
Charles Tweedie, `The Stirling Numbers and Polynomials,'' {\sl
Proceedings of the Edinburgh Mathematical Society\/ \bf 37} (1918), 2--25.

\bib
[55]\enspace
Karl Weierstrass, ``Zur Theorie den eindeutigen analytischen
Functionen,''
{\sl Mathematische Abhandlungen der Akademie der Wissenschaften zu
Berlin\/} (1876), 11--60; reprinted in his {\sl Mathematische Werke},
volume~2, 77--124. (Florian Cajori, in {\sl History of Mathematical
Notations\/ \bf 2}, cites unpublished papers of 1841 and 1859 as the
first occurrences of the notation~$\vert z\vert$; however, those
papers were not edited for publication until 1894, and they use the
notation without defining it, so their published form may differ from
Weierstrass's original.)

\bib
[56] \enspace
Christian Wiener, ``Geometrische und analytische Untersuchung der
{\it Weierstrass\/}schen Function,'' 
{\sl Journal f\"ur die reine und angewandte Mathematik\/
\bf 90} (1881), 221--252.

\vfill\eject

\def\nopagenumber{\output{\shipout\box255}}
\nopagenumber

\font\eusm=eusm10
\font\eufr=eufm10

\noindent
Note to printer: A few special symbols are used herein.

\bigskip
\halign{\qquad\qquad #\hfil\cr
{\eusm S} is uppercase script S\cr
{\eufr C} is uppercase Fraktur C\cr
{\eufr S} is uppercase Fraktur S\cr
{$\cal A$}  is uppercase script A\cr
{\eusm F} is uppercase script F\cr
{\eusm R} is uppercase script R\cr
{\eufr k} is lowercase Fraktur k\cr}

\bye